\documentclass[11pt]{amsart}
\usepackage{graphics}
\usepackage{amsmath}

\newcommand{\tc}[2]{{\left[ \begin{array}{c}#1\\#2\end{array} \right]}}

\begin{document}
\def\ZZ{\mathbb{Z}}
\def\NN{\mathbb{N}}
\theoremstyle{plain}
\newtheorem{lem}{Lemma}
\newtheorem{thm}[lem]{Theorem}
\newtheorem{prop}[lem]{Proposition}
\newtheorem{definition}[lem]{Definition}
\title{Some remarks on the characters of the general Lie superalgebra}
\author{R.C. Orellana
\and
Mike Zabrocki}
\begin{abstract}
 We compute an explicit formula the Hilbert 
(Poincar\'e) series for the ring of hook Schur functions. 
\end{abstract}
\maketitle
\centerline{\bf{Introduction}}
Berele and Regev \cite{br} defined the ring of hook Schur functions and 
showed that these functions are the characters of the general Lie superalgebra.
Since the introduction of the hook Schur functions they have been extensively
studied with respect to their combinatorial properties \cite{r1},
\cite{r2}, \cite{r3}. 
\vskip 0in
We compute the Hilbert series of this ring by giving a
generating function for the partitions which fit inside of a
$(k,\ell)$-hook.  Besides its interest as a purely combinatorial
result, this formula should have applications to discovering
algebraic properties of the ring of hook Schur functions.

\section{Definitions and Notation}
For $k$, $\ell \in \NN$ let $x_1,\ldots , x_k; y_1,\ldots , y_\ell$ be two
sets  of commuting variables.
In this note $\lambda,\  \mu,\  \nu$ will represent partitions and these
partitions will be identified with their corresponding
Young diagrams (geometric representation of partitions). 
\vskip 0in
Let $n\in \NN$. Define $H_n(k,\ell)$ to be the set of partitions
$\lambda=(\lambda_1,\cdots, \lambda_m)$ such that $\lambda_{k+i}\leq \ell$
for $i=1,\cdots, m-k$. This is equivalent to Young diagrams that fit
inside the $(k,\ell)$-hook, see Figure 1.
\vskip 0in
\begin{picture}(200,100)(0,0)
\put(150,30){\line(0,0){60}}
\put(180,30){\line(0,0){30}}
\put(150,90){\line(1,0){60}}
\put(180,60){\line(1,0){30}}
\put(210,80){\vector(0,0){10}}
\put(207,71){$k$}
\put(210,70){\vector(0,-1){10}}
\put(162,33){\vector(-1,0){12}}
\put(163,30){$\ell$}
\put(168,33){\vector(1,0){12}}
\put(160,10){Figure 1}
\end{picture}
\vskip 0in
\begin{definition}
(a)Let $\lambda$ be a partition of $n$. Then the Hook Schur function $HS_{\lambda}$ is defined as follows, for any $k$,
$\ell$,
$$
HS_{\lambda}(x_1,\ldots , x_k;y_1,\ldots, y_\ell):=
\sum_{\mu<\lambda}s_{\mu}(x_1,\ldots, x_k)s_{\lambda'/\mu'}(y_1,\cdots,
y_\ell),
$$
where $s_{\nu}$ denotes the Schur function and $\lambda'/\mu'$ denotes the 
conjugate of the skew partition $\lambda/\mu$.
\vskip .05in
(b) $\Lambda_n^{(k,\ell)}= span_F\{HS_{\lambda}(x_1,\cdots, x_n;
y_1,\ldots, y_\ell)\, |\, \lambda \in H_n(k,\ell)\}$ where $H_n(k,\ell)$
is the subset of the partitions of $n$ which fit in the 
$(k,\ell)$ hook. 
\vskip 0in
Then form $\Lambda^{(k,\ell)}= \bigoplus_{n\geq 0} \Lambda_n^{(k,\ell)}$.
\end{definition}
These functions were defined in \cite{br} as the characters of the general Lie
superalgebra, $pl(V)$, where $V$ is a vector space with an associated 
$\ZZ/2\ZZ$ grading. Berele and Regev proved that $\Lambda^{(k,\ell)}$ is
a ring.
\vskip 0.1in

\section{Hilbert Series of $\Lambda^{(k,\ell)}$}
Computing this series is equivalent to the combinatorial problem of counting 
the number of Young diagrams which fit inside a $(k,\ell)$ hook. It is
well-known that the generating function for partitions with less than or
equal to $k$  rows is given by 
$G_k(t)=\prod_{i=1}^k \frac{1}{1-t^i}$ and that the generating function
for the partitions which fit inside of a $k \times \ell$ rectangle is
given by $\tc{k+\ell}{\ell} = \frac{(t;t)_{k+\ell}}{(t;t)_k (t;t)_\ell}$
where $(a;x)_n = (1-a)(1-a x) (1-a x^2) \cdots (1-a x^{n-1})$.
\vskip 0in
Let $G_{k,\ell}(t)$ be the generating function of the partitions that 
fit in the $(k,\ell)$ hook.  Clearly we have that $G_{k,\ell}(t)$ are
symmetric with respect to the $k$ and $\ell$ parameters and it is quite
easy to see combinatorially that this function must satisfy the
recurrence
$$G_{k,\ell}(t)=G_{k,\ell-1}(t)+t^{\ell(k+1)}G_k(t)G_\ell(t)$$
since the right hand side represents the generating function of
the partitions that fit inside of the $(k,\ell)$-hook which do not contain
the cell $(k+1,\ell)$, plus the generating function for the partitions
in the hook which do contain the cell $(k+1, \ell)$.

For our computations we use mainly the following 
$t$-binomial identities.
\begin{align}
(1-t^i) \tc{k}{i} &= (1-t^{k-i+1}) \tc{k}{i-1}\\
\tc{a}{b} &= \tc{a-1}{b-1} + t^b \tc{a-1}{b} \label{idone} 
\end{align}
 
\begin{lem} \label{closed} Let $A^{k,\ell}(t)$ be defined by the
recurrence
$$A^{k,\ell}(t) =
(1-t^{k+\ell}) A^{k,\ell-1}(t) +  t^{\ell(k+1)} \tc{k+\ell}{\ell}$$
with $A^{k,0}(t) = 1$.  Then 
\begin{equation}\label{Adef}
A^{k,\ell}(t) = 1 + \sum_{i=1}^\ell t^{i(k+\ell+1)} \tc{k}{i}
\sum_{j=0}^{\ell-i} t^{j(k-i+1)} \tc{i+j-1}{j}
\end{equation}
\end{lem}

\begin{proof}
We proceed by showing that $A^{k,\ell}(t)$ given in equation
(\ref{Adef}) satisfies the relation $A^{k,\ell} - A^{k,\ell-1} +
t^{k+\ell} A^{k,\ell-1}= t^{\ell(k+1)} \tc{k+\ell}{\ell}$.  The full
details of this calculation are not especially enlightening, hence we
leave them to the reader. We give an outline of this proof by stating
that as an intermediate step for the left hand side of this equation, one
has

\begin{align*}
A^{k,\ell}(t) -& A^{k,\ell-1}(t) + t^{k+\ell}A^{k,\ell-1}(t)=  \\
&t^{\ell(k+1)} +
\sum_{i=1}^{\ell-2} t^{(i+1)(k+\ell)}
\tc{k}{i} \sum_{j=0}^{\ell-i-2} t^{(j+1)(k-i)} \tc{i+j}{j}\\
& +\sum_{i=1}^{\ell-2} t^{(i+1)(k+\ell)}\tc{k}{i}
\sum_{j=0}^{\ell-2-i} t^{j(k-i)} \left( t^j \tc{i+j-1}{j} -
\tc{i+j}{j}\right)\\ 
&+\sum_{i=1}^{\ell-2} t^{(i+1)(k+\ell)}\tc{k}{i} t^{(\ell-i-1)(k-i+1)}
\tc{\ell-2}{\ell-1-i}\\
&+ \sum_{i=1}^{\ell} t^{\ell(k+1) + i^2} \tc{k}{i} \tc{\ell-1}{\ell-i}
 + t^{\ell(k+\ell)} \tc{k}{\ell-1} \tc{\ell-2}{0}
\end{align*}

This may be reduced by repeated applications of (\ref{idone}) until one
arrives at

\begin{equation}
= t^{\ell(k+1)} \sum_{i=0}^\ell t^{i^2} \tc{k}{i} \tc{\ell}{i} =
t^{\ell(k+1)}
\tc{k+\ell}{\ell}
\end{equation}
\end{proof}

\begin{thm}
The Hilbert series of $\Lambda^{(k,\ell)}$ is given by 
$$G_{k,\ell}(t)=A^{k,\ell}(t) \prod_{i=1}^{k+\ell}\frac{1}{1-t^i}$$
\end{thm}
\begin{proof}  Take the recurrence for $G_{k,\ell}(t)$ and divide by
$G_{k+\ell}(t)$. By Lemma \ref{closed}, $A^{k,\ell}(t)$ satisfies the
desired recurrence which gives us  our result.
\end{proof}

It would be interesting to have a purely combinatorial proof of this
identity by showing that any partition that fits inside of a
$(k,\ell)$-hook is isomorphic to a partition with less than or equal to
$k+\ell$ rows and some other object counted by $A^{k,\ell}(t)$.

\bibliographystyle{plain}

\end{document}